\documentclass[12pt,reqno]{amsart}
\usepackage{amssymb,verbatim,enumerate,ifthen,hyperref}
\usepackage[mathscr]{eucal}
\usepackage[utf8]{inputenc}
\usepackage[T1]{fontenc}
\def\R{\mathbb{R}}

\def\C{\mathbb{C}}

\newtheorem{theorem}{Theorem}
\newtheorem*{theorem*}{Theorem}
\def\Thm#1#2{\ifthenelse{\equal{#1}{*}}{\begin{theorem*}#2\end{theorem*}}
             {\begin{theorem}\label{T#1}#2\end{theorem}}}
\newtheorem{Atheorem}{Theorem}

\def\thm#1{Theorem~\ref{T#1}}

\newtheorem{proposition}[theorem]{Proposition}
\newtheorem*{proposition*}{Proposition}
\def\Prp#1#2{\ifthenelse{\equal{#1}{*}}{\begin{proposition*}#2\end{proposition*}}
             {\begin{proposition}\label{P#1}#2\end{proposition}}}

\newtheorem{corollary}[theorem]{Corollary}
\newtheorem*{corollary*}{Corollary}
\def\Cor#1#2{\ifthenelse{\equal{#1}{*}}{\begin{corollary*}#2\end{corollary*}}
             {\begin{corollary}\label{C#1}#2\end{corollary}}}

\newtheorem{lemma}[theorem]{Lemma}
\newtheorem*{lemma*}{Lemma}
\def\Lem#1#2{\ifthenelse{\equal{#1}{*}}{\begin{lemma*}#2\end{lemma*}}
             {\begin{lemma}\label{L#1}#2\end{lemma}}}

\theoremstyle{definition}
\newtheorem{remark}[theorem]{Remark}
\newtheorem*{remark*}{Remark}
\def\Rem#1#2{\ifthenelse{\equal{#1}{*}}{\begin{remark*}\rm #2\end{remark*}}
             {\begin{remark}\label{R#1}\rm #2\end{remark}}}

\newtheorem{example}[theorem]{Example}
\newtheorem*{example*}{Example}
\def\Exa#1#2{\ifthenelse{\equal{#1}{*}}{\begin{example*}\rm #2\end{example*}}
             {\begin{example}\label{Ex#1}\rm #2\end{example}}}

\def\eq#1{{\rm(\ref{E#1})}}
\def\Eq#1#2{\ifthenelse{\equal{#1}{*}}
  {\begin{equation*}\begin{aligned}#2\end{aligned}\end{equation*}}
  {\begin{equation}\begin{aligned}\label{E#1}#2\end{aligned}\end{equation}}}

\begin{document}
\begin{flushright}
Publ. Math. Debrecen \textbf{81}(1-2) (2012), 243–249. \\ 
\href{http://dx.doi.org/10.5486/PMD.2012.5359}{10.5486/PMD.2012.5359}
\end{flushright}
\vspace{5mm}

\title[]{Wigner's theorem revisited}

\author[Gy. Maksa]{Gyula Maksa}
\author[Zs. P\'ales]{Zsolt P\'ales}
\address{Institute of Mathematics, University of Debrecen,
H-4010 Debrecen, Pf.\ 12, Hungary}
\email{\{maksa,pales\}@science.unideb.hu}

\subjclass[2000]{Primary 39B52, Secondary 46C99} 
\keywords{Isometry, Wigner's unitary-antiunitary theorem}

\thanks{This research has been supported by the Hungarian
Scientific Research Fund (OTKA) Grant NK81402 and by the
T\'AMOP 4.2.1./B-09/1/KONV-2010-0007 project implemented through the New
Hungary Development Plan co-financed by the European Social Fund, and the
European Regional Development Fund}

\begin{abstract}
In this paper, we give the general solution of the functional equation
$$
  \big\{\|f(x)+f(y)\|,\|f(x)-f(y)\|\big\}=\big\{\|x+y\|,\|x-y\|\big\}\qquad(x,y\in X)
$$
where $f:X\to Y$ and $X, \,\,Y$ are inner product spaces. Related equations are also considered.
Our main tool is a real version of Wigner's unitary--antiunitary theorem.
\end{abstract}

\maketitle

\section{Introduction}

An isometry from a normed space $X$ into another normed space $Y$ is a function $f:X\to Y$ which satisfies
the equality
\Eq{MU}{
  \|f(x)-f(y)\|= \|x-y\| \qquad (x,y\in X).
}
This equation implies strong structural properties for the function $f$. A classical result in this direction is
a celebrated theorem of Mazur and Ulam \cite{MazUla32} which states that an isometry $f$ of a real normed space
$\emph{onto}$ another normed space is necessarily affine. In other words, for the surjective solutions
$f:X\to Y$ of \eq{MU}, $x\mapsto f(x)-f(0)$, is a norm preserving linear map. Baker \cite{Bak71c}
showed that the same conclusion remains valid if the surjectivity assumption is replaced by the strict convexity
of the target space $Y$. Another important result which is related to linear isometries is Wigner's theorem
\cite{Wig31} and its generalization obtained by R\"atz \cite[Corollary 8(a)]{Rat96}. For further generalizations of
this fundamental result, we mention the papers \cite{AlmSha92}, \cite{Che07}, \cite{Gyo04a}, \cite{Mol96g},
\cite{Mol98a}, \cite{Mol99a}, \cite{Mol00e}, \cite{Mol00d}, \cite{Mol01c}, \cite{Mol02f}, and \cite{SimMukChaSri08}.

Assuming that $X$ and $Y$ are
\emph{real} inner product spaces, R\"atz's result characterizes functions $f:X\to Y$ that are phase equivalent
to a linear isometry (i.e., there exists a function $\varepsilon:X\to\{-1,1\}$ such that $\varepsilon f$
is a norm preserving real linear map) by the property
\Eq{2}{
|\langle f(x),f(y)\rangle|=|\langle x,y\rangle|  \qquad(x,y\in X).
}
In the complex setting, Wigner's theorem \cite{Wig31} (cf.\ also \cite{Gyo04a}) says that 
the solutions of \eq{2} are phase equivalent to a linear or conjugate linear isometry.
Without assuming that $X$ and $Y$ are real inner product spaces, we can easily see that all functions
$f:X\to Y$ that are phase equivalent to a real linear isometry are also solutions of the functional equation
\Eq{1}{
  \big\{\|f(x)+f(y)\|,\|f(x)-f(y)\|\big\}=\big\{\|x+y\|,\|x-y\|\big\}\qquad(x,y\in X).
}
Indeed, if $\varepsilon:X\to\{-1,1\}$ and $g:=\varepsilon f$ is a norm preserving real linear map,
then, for all $x,y\in X$,
\Eq{*}{
  \|f(x)\pm f(y)\|&=\|\varepsilon(x)g(x)\pm\varepsilon(y)g(y)\|
   =\|g(x)\pm\varepsilon(x)\varepsilon(y)g(y)\|\\
   &=\|g(x\pm\varepsilon(x)\varepsilon(y)y)\|=\|x\pm \varepsilon(x)\varepsilon(y)y\|,
}
which implies \eq{1} because $\varepsilon(x)\varepsilon(y)$ is either equal to $1$ or to $-1$.

The aim of this short note is to show that the converse also holds provided that $X,Y$ are inner product spaces.
That is, in that case, all solutions $f:X\to Y$ of \eq{1} are phase equivalent to a real linear isometry.
The main tool in the proof is R\"atz's characterization theorem described above.

\section{The equivalence of some functional equations related to \eq{1} and our main results}

Throughout the remaining part of this paper, $X$ and $Y$ denote real or complex inner product spaces.
We note that every complex linear space is trivially a real linear space and if $\langle\cdot,\cdot\rangle$
is a complex inner product on $X$ (or on $Y$) then $\langle\!\langle\cdot,\cdot\rangle\!\rangle$ defined
as $\langle\!\langle x,y\rangle\!\rangle=\Re\langle x,y\rangle$ is real inner product on $X$
which induces the same norm. (Here $\Re z$ stands for the real part of the complex number $z$.)
Therefore, we may assume that $\langle\cdot,\cdot\rangle$ always denotes the real inner product on $X$ and $Y$.
A function $f:X\to Y$ is called real linear if $f$ is additive and homogeneous with respect to real numbers.
Real linearity does not imply complex linearity in general as it is shown by the following example
constructed by R\"atz \cite{Rat96}: Let $X=Y=\C^{2}$ equipped with the usual inner product
\Eq{*}{
  \langle (x_{1},x_{2}),(y_{1},y_{2})\rangle=x_{1}\overline{y}_{1}+x_{2}\overline{y}_{2}
    \qquad((x_{1},x_{2}),(y_{1},y_{2})\in\C^{2}),
}
and $f(x_{1},x_{2})=(x_{1},\overline{x}_{2})$ for $(x_{1},x_{2})\in\C^{2}$.  An easy calculation shows
that $f$ is norm preserving real linear, but it is not complex-homogeneous, and hence it is not linear.

We begin with a characterization of norm-preserving real linear maps between inner product spaces.

\Thm{1}{For any $f:X\to Y$, the following three statements are equivalent:
\begin{itemize}
\item [(i)] $\|f(x)+f(y)\|=\|x+y\| \quad(x,y\in X)$;
\item [(ii)] $\langle f(x),f(y)\rangle=\langle x,y\rangle\quad(x,y\in X)$;
\item [(iii)] $f$ is a norm-preserving real linear map.
\end{itemize}}

\begin{proof}
Suppose first that (i) holds. Putting $x=y$, it follows that $f$ is norm-preserving. Now using (i)
and the norm preserving property, we get
\Eq{*}{
  2\langle f(x),f(y)\rangle
  =\|f(x)+f(y)\|^2-\|f(x)\|^2-\|f(y)\|^2=\|x+y\|^2-\|x\|^2-\|y\|^2=2\langle x,y\rangle,
}
which proves (ii).

Now suppose (ii). Putting $x=y$, the norm preserving property of $f$ follows. Using (ii) three times, 
for all $x,y,z\in X$, we obtain
\Eq{*}{
  \langle f(x+y)-f(x)-f(y), f(z)\rangle = \langle x+y, z\rangle-\langle x,z\rangle-\langle y, z\rangle=0.
}
Applying this identity for $z\in\{x+y,x,y\}$, we get
\Eq{*}{
  \langle f(x+y)-f(x)-f(y), f(x+y)-f(x)-f(y)\rangle =0,
}
which yields that $f$ is additive.

Finally, assume that $f$ is a norm-preserving real linear map. Then, by the additivity and the norm-preserving
property, we get that $\|f(x)+f(y)\|=\|f(x+y)\|=\|x+y\|$ which implies (i).
\end{proof}

\begin{remark*}
The equivalence of (i) and (iii) can easily be proved by supposing only that $X$ and $Y$ are normed
spaces and $Y$ is strictly convex. Indeed, the substitution $y=x$ in (i) implies that $f$ is norm-preserving.
Therefore $\|f(x)+f(y)\|=\|f(x+y)\|$ holds for all $x,y\in X$. Applying a result of Ger \cite{Ger93b},
we obtain that $f$ is additive which implies (iii). On the other hand, (i) follows from (iii) immediately.
\end{remark*}

In the following theorem, we list four equivalent conditions that are equivalent to \eq{1}.

\Thm{2}{For any $f:X\to Y$, the following five statements are equivalent:
\begin{itemize}
\item [(i)] \eq{1} holds;
\item [(ii)] $\|f(x)+f(y)\|+\|f(x)-f(y)\|=\|x+y\|+\|x-y\| \quad(x,y\in X)$;
\item [(iii)] $f(0)=0 \,\, \mbox{and} \,\, \|f(x)+f(y)\|\|f(x)-f(y)\|=\|x+y\|\|x-y\|\quad(x,y\in X)$;
\item [(iv)] $|\langle f(x),f(y)\rangle|=|\langle x,y\rangle|\quad(x,y\in X)$;
\item [(v)] There exists a function $\varepsilon: X\to\{-1,1\}$ such that $\varepsilon f$ is a norm-preserving real linear map.
\end{itemize}}

\begin{proof} The statement (i) implies (ii) obviously. With the substitution $y=x$, it follows from (ii) that
$f$ is norm preserving, i.e.,
\Eq{4}{
  \|f(x)\|=\|x\|\qquad(x\in X).
}
With $x=0$, this yields $f(0)=0$. Now we square the equation in (ii) to obtain
\Eq{*}{
  \|f(x)\|^2&+2\langle f(x),f(y)\rangle +\|f(y)\|^2
    + 2\|f(x)+f(y)\|\|f(x)-f(y)\|\\
    &+\|f(x)\|^2-2\langle f(x),f(y)\rangle +\|f(y)\|^2 \\
    =&\|x\|^2+2\langle x,y\rangle +\|y\|^2
    +2\|x+y\|\|x-y\|+\|x\|^2-2\langle x,y\rangle +\|y\|^2.
}
Using \eq{4}, the above equality simplifies to the second equality in (iii). Thus (ii) implies (iii).

Substituting $y=0$ into the second equation in (iii), we get \eq{4}. Squaring the  second equation in (iii) and using
\eq{4} again, we obtain that
\Eq{*}{
  \big(\|x\|^2+2\langle f(x),f(y)\rangle +\|y\|^2\big)&\big(\|x\|^2-2\langle f(x),f(y)\rangle +\|y\|^2\big)\\
  &=\big(\|x\|^2+2\langle x,y\rangle +\|y\|^2\big)\big(\|x\|^2-2\langle x,y\rangle +\|y\|^2\big).
}
This simplifies to
\Eq{*}{
  \big(\langle f(x),f(y)\rangle\big)^2=\big(\langle x,y\rangle\big)^2\qquad(x,y\in X),
}
which is equivalent to the equation in (iv) proving that (iii) implies (iv).

If (iv) holds, then, by the result of R\"atz \cite[Corollary 8(a)]{Rat96} described in the introduction, (v) follows.

Finally, (v) implies (i) as we have seen it in the introduction.
\end{proof}

The following corollary describes the continuous solutions of \eq{1}.

\Cor{5}{Let $X$ be at least two dimensional. For a continuous function $f:X\to Y$, the four equivalent
statements (i)--(iv) of \thm{2} hold if and only if $f$ is a norm-preserving real linear map.}

\begin{proof}
Assume that $f$ is a continuous function satisfying any of the conditions (i)--(iv) of \thm{2}. Then
there exists a function $\varepsilon:X\to \{-1,1\}$ such that $\varepsilon f$ is norm-preserving and real linear.
Thus, by \thm{1},
\Eq{*}{
  \varepsilon(x)\varepsilon(y)\langle f(x),f(y)\rangle=\langle x,y\rangle\qquad(x,y\in X).
}
If $y\neq0$, then there exists an open ball $U$ around $y$ such that
\Eq{*}{
  \varepsilon(x)=\varepsilon(y)\frac{\langle x,y\rangle}{\langle f(x),f(y)\rangle}\qquad(x\in U).
}
This, by the continuity of $f$, shows that $\varepsilon$ is continuous on $U$ and hence it is constant on $U$.
The set $X\setminus \{0\}$ is connected (because $X$ is at least two dimensional), therefore $\varepsilon$
is constant on $X\setminus \{0\}$. Thus $f$ must be a norm-preserving real linear map.
\end{proof}

\begin{remark*}
In the exceptional nontrivial case when $X$ is one dimensional and real, say $X=\{\lambda a: \, \lambda\in\R\}$ with
some $a\in X,\, \|a\|=1$, and $Y$ is at least one dimensional, the above argument shows that $\varepsilon$ is constant on
the set of positive reals and constant also on the set of negative reals. Therefore $f$ is either a norm-preserving real
linear map or $f(\lambda a)=|\lambda|b$ for all $\lambda\in\R$ and for some $b\in Y$ with $\|b\|=1$.
\end{remark*}

Finally, we formulate two open problems.

\textbf{Problem 1.} What are the solutions $f:X\to Y$ of \eq{1} when $X$ and $Y$ are normed but not
necessarily inner product spaces? Under what conditions does it remain valid that, for the solutions
of \eq{1}, $\varepsilon f$ is real linear for some function $\varepsilon:X\to\{-1,1\}$?

\textbf{Problem 2.} Let $X$ and $Y$ be complex normed spaces. Let $n$ be a fixed positive integer and denote
$\beta_{1}, \dots, \beta_{n}$ the $n$th roots of unity. These elements form a multiplicative subgroup of the 
unit circle in $\C$. Find the solutions $f:X\to Y$ of the following generalization of \eq{1}:
\Eq{22}{
\big\{\|f(x)-\beta_{k}f(y)\|: k\in\{1,\dots, n\}\big\}=\big\{\|x-\beta_{k}y\|: k\in\{1,\dots, n\}\big\}  \qquad(x,y\in X).
}
Obviously, this is the isometry equation in case $n=1$, and the case $n=2$ was just discussed in this paper.
One can also see that if there exists a function $\varepsilon:X\to \{\beta_{1}, \dots, \beta_{n}\}$ such that
$\varepsilon f$ is complex linear and norm-preserving, then $f$ satisfies \eq{22}. Under what conditions 
does it remain valid that, for the solutions of \eq{22}, $\varepsilon f$ is complex linear and norm-preserving 
for some function $\varepsilon:X\to \{\beta_{1}, \dots, \beta_{n}\}$?

\bigskip

\textbf{Acknowledgement.} The authors are indebted to Professor Lajos Molnár, who 
suggested the investigation of the functional equation in (ii) of Theorem 2 and had also 
several valuable comments.


\end{document}